\documentclass[a4paper,10pt]{article}
\usepackage[pdftex]{graphicx}
\usepackage{setspace}
\usepackage{amssymb, amsthm, amsmath, enumerate, graphicx, epsfig, hyperref,geometry,color}
\usepackage[T1]{fontenc}
\def\name{Junnan He}
\def\title{\huge{An alternative definition for the Conley relation}}
\hypersetup{
  colorlinks = true,
  urlcolor = cyan,
  pdfborder={text}
  pdftoolbar = true,
  pdftitle={\title},
  pdfauthor = {\name},
  pdfkeywords = {Cauchy-Schwarz, A.M.-A.M., H$\ddot{\text{o}}$lder, Convexity, Mid-Convexity, Theory of Convexity, Cesaro Mean, Stolz-Cesaro Theorem, },
  pdftitle = {\name: An alternative definition for the Conley relation},
  pdfsubject = {Relativistic Jet Modeling},
  pdfpagemode = UseNone,
  citecolor=blue,
  linkcolor= red,
  linkbordercolor={1 0 0}
}
\theoremstyle{definition}

\begin{document}
\fontsize{13}{16}\selectfont
\begin{center}
{\Large \title}
\ \\
\ \\
\ \\
{\large \textit \name}\footnote{Mathematical Sciences Institute, the Australian National University, Canberra ACT Australia. Contact: \href{.../}{\tt u4575078@anu.edu.au}}\ \\
\end{center}
\ \\
\ \\
\begin{abstract}
In the theory of dynamics of closed relations on compact Hausdorff spaces, the definition for the Conley relation $f^\Omega$ of an iterated closed relation $f$ is nontrivial. This paper establishes a new equivalent definition for $f^\Omega$, and discusses its interpretation.
\end{abstract}\

\section{Introduction}\
\\
The Conley Decomposition Theorem has been referred to as the fundamental theorem of dynamical systems \cite{Norton}. Originally it describes the asymptotic behavior of a dynamical system on a compact metric space. In constructing the proof, Conley \cite{Conley} defined a weak form of recurrence coined as ``chain recurrence''. The idea was that a point is chain recurrent if it can always return to itself via some long orbits so long as making arbitrarily small errors along the way is allowed.

It has been argued lately that the natural settings for the Conley Decomposition Theorem to study the decomposition of a compact Hausdorff space under iterations of closed relations (\cite{Akin},  \cite{McGehee}). Based on the original development by Conley \cite{Conley}, McGehee and Wiendt \cite{McGehee} naturally generalized the notion of Conley relation to the study of iterated relations on compact topological spaces. The Conley relation $f^\Omega$ of a relation $f$ can be thought of as the ``infinite iterate'' of $f$ that takes into account a topological structure on the space of relations \cite{McGehee}.  In developing $f^\Omega$, two other relations $f^\infty$ and $f^\omega$ were introduced. The three relations can all be thought of as ``infinite iterates'' of $f$, but differences lie in that $f^\infty$ takes no topology into account and $f^\omega$ takes into account only the topology of the phase space. Under the construction, a point is chain recurrent in the iterated relations sense if and only if it is related to itself by the Conley relation.

The extension given by McGehee and Wiendt \cite{McGehee} was very general, but the definition for the extended Conley relation seems more sophisticated than the original one because of the topological structure on the space of relations. In this paper, we observe some interesting interplay of the relations $f^\Omega$ and $f^\omega$, and we give an interpretation of the Conley relation by proving an alternative definition for $f^\Omega$. Many properties of the three limiting relations were proven by McGehee and Wiendt \cite{McGehee}. Our proof would be heavily dependent upon those properties as well as some other topological results about iterated relations.

In this paper, we first explain the motivations and intuitions of the alternative definition. Then we state some of the useful results given by McGehee and Wiendt \cite{McGehee}. In Section 4, we would used these results to prove the alternative definition for $f^\Omega$.

\section{Notations and Motivations}\
\\
Consider $X$ as a compact Hausdorff topological space and $f \subset X^2$ is a closed relation on $X$. McGehee and Wiandt \cite{McGehee} studied the dynamics of $X$ under the iteration of $f$ in detail. They also adapted and proposed some definitions of the limiting relations that are very useful in studying the dynamics. In this paper, we would adopt the notations used in McGehee and Wiandt. For a subset $A$ of a topological space $X$, a (open/closed) neighbourhood $U$ of $A$ is a (open/closed) set that contains an open set containing $A$. We write the family of open neighbourhoods of $A$ as $\mathfrak N^{o} (A)$, the family of neighbourhoods as $\mathfrak N  (A)$ and the family of closed neighbourhoods as $\overline{\mathfrak N } (A)$. The following definitions are given by McGehee and Wiandt \cite{McGehee}. \\
\\
{\bf{Definition\ 2.1} } If $f$ is a relation on a set $X$, then the \emph{limit relation} of $f$ is
\[ f^{\infty} \equiv \bigcap_{n \geq 0} \bigcup_{k\geq n} f^k \]
\\
{\bf{Definition\ 2.2} } If $f$ is a relation on a set $X$, then the \emph{$\omega$-limit relation} of $f$ is
\[ f^{\omega} \equiv \bigcap_{n \geq 0} \overline{\bigcup_{k\geq n} f^k} \]
\\
{\bf{Definition\ 2.3} } If $f$ is a relation on a set $X$, then the \emph{Conley relation} of $f$ is
\[ f^{\Omega} \equiv \bigcap_{\phi \in \overline{\mathfrak N}(f)} \phi^\infty \]

Among the three limiting relations, the Conley relation is the most fundamental since it is closedly related to the chain recurrent sets. For example, it is shown by McGehee and Wiandt \cite{McGehee} that if $f$ is a closed relation on a compact Hausdorff space $X$, then the chain recurrent set of $X$ under $f$ is the set of fixed points under $f^{\Omega}$. Intuitively, the above definition for the Conley relation should have been modified from the conventional concepts of recurrence in the literature of dynamical systems. And conventionally, the chain recurrence on a metric space is defined using $\epsilon$-pseudo-orbits. The following definition was adopted by Norton \cite{Norton}.\\
\\
{\bf{Definition\ 2.4}} Given $x, y \in X$ and $\epsilon >0$, an \emph{$\epsilon$-pseudo-orbit} from x to y means a sequence of points $(x = x_0, x_1 , \dots , x_n =y)$ with $n>0$ such that for $k=0,1, \dots , n-1$, we have $d(f(x_k), x_{k+1}) < \epsilon$.

\begin{center}\includegraphics[width=140mm]{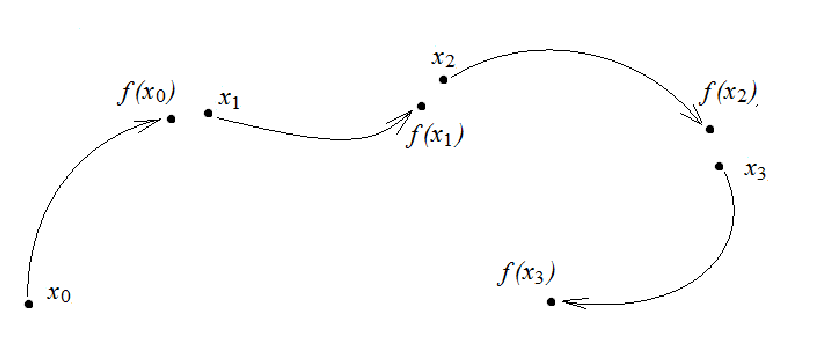} \\
\begin{flushleft}
\small{$\bf{Figure\ 1.} $ This picture is a simple illustration of an $\epsilon$-pseudo-orbit from the point $x_0$ to $x_3$.}
\end{flushleft}
\end{center}

In light of Definition 2.4 and Figure 1, it seems to suggest that to find points in an $\epsilon$-pseudo-orbit of $x \in X$, we need to think of the dynamical system under the iteration of some $\iota_{\epsilon} \circ f$ where $\iota_{\epsilon}$ maps each point of $X$ into its $\epsilon$-neighbourhood. Such $\iota_{\epsilon} \circ f$ should be a subset of some $\phi \in \mathfrak N(f)$. It can be shown that in a metric space, for every $\iota_{\epsilon}$, there exists $\phi \in \mathfrak N(f)$ such that $ \iota_{\epsilon} \circ f \subset \phi$ in terms of the product topology in $X^2$, and if $f$ is closed and $X$ is compact, every $\phi \in \mathfrak N(f)$ contains some $\iota_{\epsilon} \circ f$. This gives an interpretation for Definition 2.3 that if for every $\epsilon >0$, there is an  $\epsilon$-pseudo-orbit from $x$ to $y$, then $y\in f^{\Omega}(x)$. But the interpretation seems to suggest more. For any $x\in X$, the smaller $\epsilon$ we take, the closer the $\epsilon$-pseudo-orbit of $x$ would line up with the actual orbit. Because the actual orbit approaches $f^{\omega}(x)$, chances are the $\epsilon$-pseudo-orbit would meet $f^{\omega}(x)$ as well. If for any neighbourhood of $f^{\omega}(x)$, there is an $\epsilon$-pseudo-orbit that meets it, then heuristically we can calculate $f^{\Omega}(x)$ by looking for points lying in pseudo-orbits starting from $f^{\omega}(x)$.

In this paper, we would use $I \subset X^2$ to denote the diagonal elements, i.e. $I \equiv \{(x,x)| x\in X\}$. Standard topology argument shows if $X$ is compact, $I$ is a closed subset of $X^2$. In addition, we would use $\iota$ to denote a neighbourhood of $I$, i.e. $\iota \in \mathfrak N(I)$.

\section{Some Useful Results}\
\\
The following constructions are given in McGehee and Wiandt \cite{McGehee}, and will be referred to frequently in the next section.\\
\\
{\bf{Definition\ 3.1}} If $\mathfrak K$ is a family of subsets of a fixed set, then $\mathfrak K$ is called \emph{directed} if and only if the following statement holds. If $K_1 \in \mathfrak K$ and $ K_2 \in \mathfrak K$, then there exists $K_3 \in \mathfrak K$ satisfying $K_3 \subset K_1 \cap K_2$.\\
\\
{\bf{Definition\ 3.2}} If $X$ and $Y$ are topological spaces, then a map $\Theta: 2^X \times 2^X \rightarrow 2^Y$ is called \emph{semicontinuous} if and only if it satisfies the following conditions.\\
(a) $\Theta$ preserves inclusion.\\
(b) $\Theta$ preserves closed sets.\\
(c) If $\mathfrak K$ and $\mathfrak L$ are directed families of closed subsets of $X$, then
\[ \Theta \left(\bigcap \mathfrak K, \bigcap \mathfrak L \right) = \bigcap \{\Theta(K,L)| K\in \mathfrak K, L \in \mathfrak L\}. \]
\\
{\bf{Definition\ 3.3}} Let $\mathcal R(X)$ be the family of relations on the topological space $X$. The map
\[ \mathcal R(X) \times \mathcal R(X) \rightarrow \mathcal R(X): (f,g) \rightarrow g \circ f\]
will be called the \emph{composition map}.\\
\\
{\bf{Theorem\ 3.4}} (Theorem 4.8 in \cite{McGehee}) If $X$ is a compact Hausdorff space, then the composition map is semicontinuous.\\
\\
{\bf{Theorem\ 3.5}} (Theorem 8.5 in \cite{McGehee}) If $f$ is a closed relation on a compact Hausdorff space $X$ , then the following properties hold.\\
(a) $f \circ f^\Omega = f^\Omega$.\\
(b) $f^\Omega \circ f = f^\Omega$.\\
(c) $f^\Omega \circ f^\Omega = f^\Omega$.\\
\\
{\bf{Theorem\ 3.6}} (Theorem 8.3 in \cite{McGehee}) If $f$ is a closed relation on a compact Hausdorff space $X$, then the following properties hold.\\
(a) $f^\omega \subset f^\Omega.$\\
(b) $ f^\Omega = \bigcap \{\phi^\Omega | \phi \in \overline{\mathfrak N}(f)\}.$ \\
(c) $ f^\Omega = \bigcap \{\phi^\omega | \phi \in \overline{\mathfrak N}(f)\}.$ \\
(d) $f^\Omega$ is closed.\\
\\
{\bf{Lemma\ 3.7}} (Lemma 7.5 in \cite{McGehee}) If $f$ is a closed relation on a compact Hausdorff space $X$ and $g\in \mathfrak N(f)$, then $g^\infty \in \mathfrak N(f^\omega).$\\
\\
{\bf{Lemma\ 3.8}} (Lemma 5.5 in \cite{McGehee}) If $f$ is a closed relation on a compact Hausdorff space $X$ and $g\in \mathfrak N(f)$, then there exist $\iota', \iota'' \in \mathfrak N(I)$ such that $\iota' \circ f \circ \iota'' \subset g$.\\
\\
{\bf{Theorem\ 3.9}} (Theorem 7.3 in \cite{McGehee}) If $f$ is a closed relation on a compact Hausdorff space $X$, then the following inclusions hold.\\
(a) $f \circ f^\omega \supset f^\omega$.  \\
(b) $f^\omega \circ f \supset f^\omega$.  \\
(c)  $f^\omega \circ f^\omega \supset f^\omega$.\\
\\
{\bf{Lemma\ 3.10}} (Lemma 4.13 (b) in \cite{McGehee}) If $f, g$ and $h$ are relations on a topological space $X$, then for any $\phi \in \mathfrak N(f)$, $\psi \in \mathfrak N(h)$, $\psi \circ g \circ \phi \in \mathfrak N(h\circ g \circ f)$.\\

\section{The Main Result}
{\bf{Lemma\ 4.1}} If $f$ is a closed relation on a compact Hausdorff space $X$, then $f^{\Omega} \circ f^{\omega} = f^{\Omega}$.
\begin{proof} \begin{align*}
f^{\Omega} \circ f^{\omega} &= f^\Omega \circ \bigcap_{n \geq 0} \overline{\bigcup_{k\geq n} f^k} = \bigcap_{n \geq 0} \left( f^\Omega \circ \overline{\bigcup_{k\geq n} f^k} \right) & \text{by Theorem 3.4}\\
& \supset \bigcap_{n \geq 0} \left( f^\Omega \circ \bigcup_{k\geq n} f^k \right)  =\bigcap_{n \geq 0}   \bigcup_{k\geq n} \left( f^\Omega \circ f^k \right)  & \text{by Theorem 3.4}\\
& =\bigcap_{n \geq 0}   \bigcup_{k\geq n} \left( f^\Omega  \right) & \text{by repeatedly apply Theorem 3.5 (b)}\\
&= f^\Omega  =  f^\Omega \circ f^\Omega & \text{by Theorem 3.5 (c)}\\
& \supset f^\Omega \circ f^\omega. & \text{by Theorem 3.6 (a)}\\
\end{align*}
Therefore, $f^\Omega \circ f^\omega = f^\Omega$ as desired.
\end{proof}\
\\
{\bf{Theorem\ 4.2}} (Main result) If $f$ is a closed relation on a compact Hausdorff space $X$, then the following properties hold.\\
(a) $ f^\Omega = \bigcap_{\phi \in \mathfrak N (f)} \bigcup_{n\geq 0} \left( \phi^n \circ f^\omega\right).$ \\
(b) $ f^\Omega =  \bigcap_{\iota \in \mathfrak N (I)} \bigcup_{n\geq 0} \left( (\iota \circ f)^n \circ f^\omega \right).$ \\

\begin{proof}
We establish (a) by first proving $\bigcap_{\phi \in \mathfrak N}(f) \bigcup_{n\geq 0} \left(\phi^n \circ f^\omega \right) \supset f^\Omega$. Standard topological argument and Theorem 3.4 show that 
\begin{align*}
\bigcap_{\phi \in \mathfrak N (f)} \bigcup_{n\geq 0} \left( \phi^n \circ f^\omega\right) &=\bigcap_{\phi \in \mathfrak N^o (f)} \bigcup_{n\geq 0} \left( \phi^n \circ f^\omega\right)\\
&= \bigcap_{\phi \in \mathfrak N^o (f)} \left( (\bigcup_{n\geq 0} \phi^n ) \circ f^\omega \right) & \text{by Theorem 3.4}\\
& \supset \bigcap_{\phi \in \mathfrak N^o (f)} \left( \phi^{\infty} \circ f^\omega \right)  & \text{since $ \bigcup_{n\geq 0} \phi^n \supset \phi^{\infty}$}.
\end{align*}
Because $X$ is compact Hausdorff, $X^2$ is also compact Hausdorff. Hence standard argument gives for each $\phi \in \mathfrak N^o (f), \exists \psi \in \overline{\mathfrak N} (f)$ such that $\phi \supset \psi$. Since $\phi \in \mathfrak N(\psi)$ where $\psi$ is closed, Lemma 3.7 implies $\phi^\infty \supset \psi^\omega$. Hence for each $\phi \in \mathfrak N^o (f)$, there exists $ \psi \in \overline{\mathfrak N} (f)$ such that  $\phi^\infty \supset \psi^\omega$. Therefore, 

\[ \bigcap_{\phi \in \mathfrak N^o (f)} \left( \phi^{\infty} \circ f^\omega \right) \supset  \bigcap_{\psi \in \overline{\mathfrak N }(f)} \left( \psi^{\omega} \circ f^\omega \right) . \]

In addition, since $X^2$ is compact Hausdorff, for every $\psi_1 , \psi_2 \in \overline{\mathfrak N } (f)$ there exists $\psi_3 \in \overline{\mathfrak N } (f)$ such that $\psi_1 \cap \psi_2 \supset \psi_3$ and hence $\psi_1^\omega \cap \psi_2^\omega \supset \psi_3^\omega$. So $ \{\psi^{\omega}| \psi \in \overline{\mathfrak N } (f) \} $ is a directed family of closed sets and Theorem 3.4 implies
\begin{align*}
\bigcap_{\psi \in \overline{\mathfrak N }(f)} \left( \psi^{\omega} \circ f^\omega \right) & = \left(    \bigcap_{\psi \in \overline{\mathfrak N }(f)}  \psi^{\omega}       \right) \circ f^\omega \\
& = f^\Omega \circ f^\omega = f^\Omega. & \text{by Theorem 3.6 and Lemma 4.1}
\end{align*}Therefore we have $\bigcap_{\phi \in \overline{\mathfrak N}(f)} \bigcup_{n\geq 0} \left(\phi^n \circ f^\omega \right) \supset f^\Omega$.\\

Next we want to show $f^\Omega \supset \bigcap_{\phi \in \overline{\mathfrak N}(f)} \bigcup_{n\geq 0} \left(\phi^n \circ f^\omega \right)$. By Theorem 3.6 and repeatedly applying Theorem 3.5 (a), we have for every $n \in \mathbb N$,
\[f^\Omega   = \bigcap_{\phi \in \overline{\mathfrak N}(f)}\phi^\Omega   = \bigcap_{\phi \in \overline{\mathfrak N}(f)} \left(\phi^n \circ \phi^\Omega \right)= \bigcap_{\phi \in \overline{\mathfrak N}(f)} \bigcup_{n\geq 0} \left(\phi^n \circ \phi^\Omega \right).  \]
Since $\phi \supset f$, Theorem 3.6 implies $\phi^\Omega = \bigcap_{\psi \in \overline{\mathfrak N}(\phi)}\psi^\Omega \supset  \bigcap_{\psi \in \overline{\mathfrak N}(f)}\psi^\Omega = f^\Omega$, Theorem 3.4 and Theorem 3.6 together give
\[ \bigcap_{\phi \in \overline{\mathfrak N}(f)} \bigcup_{n\geq 0} \left(\phi^n \circ \phi^\Omega \right)\supset  \bigcap_{\phi \in \overline{\mathfrak N}(f)} \bigcup_{n\geq 0} \left(\phi^n \circ f^\Omega \right) \supset  \bigcap_{\phi \in \overline{\mathfrak N}(f)} \bigcup_{n\geq 0} \left(\phi^n \circ f^\omega \right)\]
Therefore, we have shown that $f^\Omega \supset \bigcap_{\phi \in \overline{\mathfrak N}(f)} \bigcup_{n\geq 0} \left(\phi^n \circ f^\omega \right).$ This establishes (a).\\

To prove (b), it suffices to show (a) implies (b). Lemma 3.8 gives for every $\phi \in \mathfrak N(f), $ there exists $ \iota', \iota'' \in \mathfrak N(I)$ such that
\[ \phi \supset \iota' \circ f \circ \iota'' \supset \iota' \circ f \circ I = \iota' \circ f .\]
Hence
\[  \bigcap_{\phi \in \mathfrak N (f)} \bigcup_{n \geq 0} (\phi^n \circ f^\omega) \supset  \bigcap_{\iota \in \mathfrak N (I)}\bigcup_{n \geq 0} \left((\iota \circ f)^n \circ f^\omega \right).\]

We are now left to show the other inclusion, namely $\bigcap_{\phi \in \mathfrak N (f)} \bigcup_{n \geq 0} (\phi^n \circ f^\omega) \subset   \bigcap_{\iota \in \mathfrak N (I)}\bigcup_{n \geq 0} \left((\iota \circ f)^n \circ f^\omega \right).$\\

Lemma 3.8 gives that for any $\iota \in \mathfrak N(I)$, there exists $ \iota_1, \iota_2 \in \mathfrak N(I)$ such that $\iota_1 \circ I \circ \iota_2 \subset \iota$. Take $\iota' = \iota_1 \cap \iota_2$. Since $I$ is closed in compact Hausdorff $X^2$, $\iota \in \mathfrak N(I)$ and $\iota \supset \iota_1\circ \iota_2 \supset \iota' \circ \iota'$. Therefore
\begin{align*}
 \bigcap_{\iota \in \mathfrak N (I)}\bigcup_{n \geq 0} \left((\iota \circ f)^n \circ f^\omega \right) &\supset  \bigcap_{\iota' \in \mathfrak N (I)}\bigcup_{n \geq 0} \left((\iota' \circ \iota' \circ f)^n \circ f^\omega \right)\\
&=   \bigcap_{\iota' \in \mathfrak N (I)}  \left[ f^\omega \cup \bigcup_{n \geq 1} \left( \iota' \circ (\iota'  \circ f \circ \iota')^{n-1} \circ \iota' \circ f \circ f^\omega \right) \right] \\
& \supset   \bigcap_{\iota' \in \mathfrak N (I)} \left[   f^\omega \cup \bigcup_{n \geq 1} \left(  (\iota'  \circ f \circ \iota')^{n-1} \circ  f \circ f^\omega \right) \right]\\
&  \supset   \bigcap_{\iota' \in \mathfrak N (I)} \left[   f^\omega \cup \bigcup_{n \geq 1} \left(  (\iota'  \circ f \circ \iota')^{n-1}   \circ f^\omega \right) \right]
\end{align*}
by Theorem 3.9 (a). Lemma 3.10 implies that $\iota'  \circ f \circ \iota'$ is a neighbourhood of $f$. Because $X^2$ is also compact Hausdorff and $f$ is closed, there exists $\phi \in \mathfrak N(f)$ such that $\phi \subset \iota'  \circ f \circ \iota'$. It follows that
\begin{align*}
\bigcap_{\iota' \in \mathfrak N (I)} \left[   f^\omega \cup \bigcup_{n \geq 1} \left(  (\iota'  \circ f \circ \iota')^{n-1}   \circ f^\omega \right) \right] &\supset \bigcap_{\phi \in \mathfrak N (f)} \left[   f^\omega \cup \bigcup_{n \geq 1} \left( \phi^{n-1}   \circ f^\omega \right) \right]\\
& =  \bigcap_{\phi \in \mathfrak N (f)}   \bigcup_{n \geq 0} \left( \phi^{n}   \circ f^\omega \right).
\end{align*}
This gives $ \bigcap_{\iota \in \mathfrak N (I)}\bigcup_{n \geq 0} \left((\iota \circ f)^n \circ f^\omega \right) \supset \bigcap_{\phi \in \mathfrak N (f)} \bigcup_{n \geq 0} (\phi^n \circ f^\omega).$ And so we conclude that 
\[  f^\Omega = \bigcap_{\phi \in \mathfrak N (f)} \bigcup_{n\geq 0} \left( \phi^n \circ f^\omega\right) =\bigcap_{\iota \in \mathfrak N (I)} \bigcup_{n\geq 0} \left( (\iota \circ f)^n \circ f^\omega \right).\]
\end{proof}

According to the above theorem, the Conley relation is a set-valued map that first maps every point of the space to its $f^\omega$ limit, then to the set of points on ``$\epsilon$-pseudo-orbits'' from its $f^\omega$ limit. This definition gives a different angle to look at the limiting behavior in a dynamical system, and it also suggests a way in computing limits and attractors for given compact Hausdorff dynamical systems. To compute the $f^\Omega$ limit of a set in the space, we can just to calculate the  $f^\omega$ limit which may be easily found, and then consider the ``$\epsilon$-pseudo-orbits'' from  the  $f^\omega$ limit.

\section{Acknowledgements} The author wish to thank Michael F. Barnsley for all the suggestions he made in writing this report, and the wonderful lectures he taught that initiated my understandings about dynamical systems.

\end{document}